\documentclass[12pt]{article}

\usepackage{amsmath, epsfig, cite,colordvi}
\usepackage{amssymb}
\usepackage{amsfonts}
\usepackage{latexsym}
\usepackage{graphicx}

\newtheorem{thm}{Theorem}[section]

\newtheorem{cor}[thm]{Corollary}

\newtheorem{lem}[thm]{Lemma}

\newcommand{\qed}{{\hfill\rule{4pt}{7pt}}}
\def\pf{\noindent {\it Proof.} }

\def\cp{_c\partial}

\numberwithin{equation}{section}

\makeatletter \@addtoreset{equation}{section} \makeatother

\begin{document}
\begin{center}
{\bf \large An $4n$-point Interpolation Formula for \\
Certain Polynomials}\\
Sandy  H.L. Chen$^1$, Amy M. Fu$^2$\\
Center for Combinatorics, LPMC-TJKLC\\
Nankai University, Tianjin 300071, P.R. China\\ \ \\
 $^1$chenhuanlin@mail.nankai.edu.cn,
$^2$fu@nankai.edu.cn
\end{center}

\noindent{\bf Abstract}  By using some techniques of the divided difference operators,  we establish an $4n$-point interpolation formula. Certain polynomials, such as Jackson's $_8\phi_7$  terminating summation formula, are special cases of this formula.
Based on Krattenthaler's identity,  we  also give Jackson's formula  a determinantal interpretation.

\section{Introduction}

Recall that the $i$-th divided difference operator $\partial_i$, acting on functions $f(x_1,\ldots,x_n)$ of several variables, is defined by
$$
f(x_1,\ldots,x_i,x_{i+1},\ldots){\partial_i}=\frac{f(x_1,\ldots,x_{i},x_{i+1},\ldots)-f(x_1,\ldots,x_{i+1},x_{i},\ldots)}
{(x_{i}-x_{i+1})}.
$$

To be more general,  we introduce the operator  $\cp_i$ which we  call the {\it $i$-th $c$-divided difference operator}:
\begin{equation}\label{cddo}
f(x_1,\ldots,x_i,x_{i+1},\ldots)\cp_i=\frac{f(x_1,\ldots,x_{i},x_{i+1},\ldots)-f(x_1,\ldots,x_{i+1},x_{i},\ldots)}
{(x_{i}-x_{i+1})(1-c/x_ix_{i+1})}.
\end{equation}
Note that $_0\partial_i=\partial_i$.

Several  properties of  the $c$-divided difference operators will be given in next section, see Lemmas  2.1 to 2.5.
Employing those lemmas,  we obtain the main result of this paper.
\begin{thm}\label{maint}
Given two sets of variables
$$ A=\{a_1,c/a_1,\ldots,a_n,c/a_n\},\quad B=\{b_1,c/b_1,\ldots,b_n,c/b_n\},
$$ we have the following $4n$-point interpolation
formula for certain polynomials $f(y)$ of degree $2n$ with symmetry $y^{-n}f(y)=(c/y)^{-n}f(c/y)$ when $c\neq 0$:
\begin{eqnarray}
f(y)&=&\frac{f(b_1)}{\prod_{i=1}^n(b_1-a_i)(b_1-c/a_i)}\prod_{i=1}^n(y-a_i)(y-c/a_i)\nonumber\\
&&+
\frac{f(a_1)}{\prod_{i=1}^n(a_1-b_i)(a_1-c/b_i)}\prod_{i=1}^n(y-b_i)(y-c/b_i)\nonumber\\
&&+\sum_{j=1}^{n-1}C_j
\cdot \prod_{i=1}^j(y-b_i)(y-c/b_i)\prod_{i=1}^{n-j}(y-a_i)(y-c/a_i),
\end{eqnarray}
where
$$
C_j=\frac{f(b_1)b_1^{1-j}}{\prod_{i=1}^{n-j+1}(b_1-a_i)(b_1-c/a_i)} {\cp_1}\cdots {\cp_j}
(b_{j+1}-a_{n-j+1})(1-c/a_{n-j+1}b_{j+1}).
$$
\end{thm}

Note that Theorem \ref{maint} leads to the
following $2n$-point interpolation formula given in \cite{chenfu} if $c=0$:
\begin{multline*}
f(y)=\frac{f(b_1)}{\prod_{i=1}^n(b_1-a_i)}\prod_{i=1}^n(y-a_i)
+\frac{f(a_1)}{\prod_{i=1}^n(a_1-b_i)}\prod_{i=1}^n(y-b_i)\\
+\sum_{j=1}^{n-1}\frac{f(b_1)}{\prod_{i=1}^{n-j+1}(b_1-a_i)} {\partial_1}\cdots {\partial_j}
(b_{j+1}-a_{n-j+1})
\cdot \prod_{i=1}^j(y-b_i)\prod_{i=1}^{n-j}(y-a_i).
\end{multline*}

The symmetry $y^{-n}f(y)=(c/y)^{-n}f(c/y)$  when $c\neq 0$ implies that $f(y)$ can be written as a product $\prod_{i=1}^n(y-x_i)(c-x_iy)$.
Considering the case
$$
A=\{a,c/a,\ldots,aq^{1-n},cq^{n-1}/a\}, \quad B=\{b,c/b,\ldots,bq^{1-n},cq^{n-1}/b\}
$$
and $y=x_{n+1}$,  one can check that Theorem \ref{maint} implies the following identity by expanding the determinant with respect to the last row:
\begin{multline}\label{Kara}
\det\left(P_{n-j+1}(x_i,aq^{j-n})P_{n-j+1}(x_i,c/a)P_{j-1}(x_i,bq^{1-n})P_{j-1}(x_i,cq^{n-j+1}/b)\right)_{i,j=1}^{n+1}\\
=\prod_{1\leq i<j \leq n+1}(x_i-x_j)(c-x_ix_j)b^{{n+1 \choose
2}}q^{-(n+1)n(n-1)/3} \prod_{i=1}^{n+1}
(a/b,cq^{2n+2-2i}/ab;q)_{i-1},
\end{multline}
where $(a;q)_n$ is the $q$-shifted factorial  defined by
$$
(a;q)_0=1, \quad (a;q)_n=(1-a)(1-aq)\cdots(1-aq^{n-1}), \quad n=1,2,\ldots
$$
and we use $P_n(a,b)$ called Cauchy polynomials \cite{chen}
 to denote $a^n(b/a;q)_n$ for convenience.

After some rearrangement of \eqref{Kara}, we may get Krattenthaler's
identity  \cite{Warnaar}:
\begin{multline*}
\det\left(\frac{(ax_i,ac/x_i;q)_{n-j}}{(bx_i,bc/x_i;q)_{n-j}}\right)_{i,j=1}^n\\
=\prod_{1\leq i<j \leq n}(x_j-x_i)(1-c/x_ix_j)a^{{n \choose
2}}q^{\binom{n}{3}} \prod_{i=1}^n
\frac{(b/a,abcq^{2n-2i};q)_{i-1}}{(bx_i,bc/x_i;q)_{n-1}}.
\end{multline*}

Through  the specializations  $f(y)=
\prod_{i=1}^n(uq^{i-1}-y)(c-uq^{i-1}y)$ and
$$
A=\{a,c/a,\ldots,aq^{1-n},cq^{n-1}/a\}, \quad B=\{b,c/b,\ldots,bq^{1-n},cq^{n-1}/b\}
$$
in Theorem \ref{maint},  we have the following variation of Jackson's $_8\phi_7$ terminating formula.
\begin{cor} \label{coro1}   We have
\begin{align}
&\prod_{i=1}^n (uq^{i-1}-y)(c-uq^{i-1}y) \nonumber\\
&=q^{{n \choose 2}}\sum_{k=0}^n {n \brack k}\frac{P_{n-k}(b,uq^{n-1})P_k(a,uq^{n-k})P_{n-k}(u,c/b)P_{k}(c/a,u)}{P_n(b,a)(cq^{n-k-1}/ab;q)_{n-k}(cq^{2n-2k}/ab;q)_k} \nonumber\\
&\quad \times P_{n-k}(y,aq^{k+1-n})P_{n-k}(y,c/a)P_{k}(y,bq^{1-n})P_{k}(y,cq^{n-k}/b),
\end{align}
where ${n \brack k}$ is the $q$-binomial coefficient defined by
$$
{n \brack k}=\frac{(q;q)_n}{(q;q)_k(q;q)_{n-k}}.
$$
\end{cor}


In section 3,  we shall  give Corollary \ref{coro1} a different approach by considering a special case of \eqref{Kara}.  The key step of our approach is to evaluate the cofactor of each entry in the last row of the determinant in \eqref{Kara}.

Note that if we write $cq^{-1}/ab$ as $a$, $y/a$ as $b$, $c/bu$ as $c$, $c/ay$ as
$d$ and $uq^{n-1}/b$ as $e$ in the above corollary, we get
Jackson's $_8\phi_7$ formula\cite{GasperRahman}:
\begin{multline*}
\frac{(aq,aq/bc,aq/bd,aq/cd;q)_n}{(aq/b,aq/c,aq/d,aq/bcd;q)_n}
\\=\sum_{k=0}^n
\frac{(1-aq^{2k})(a;q)_k(b;q)_k(c;q)_k(d;q)_k(e;q)_k(q^{-n};q)_kq^k}{(1-a)(q;q)_k(aq/b;q)_k(aq/c;q)_k(aq/d;q)_k(aq/e;q)_k(aq^{n+1};q)_k}.
\end{multline*}

\section{The $4n$-point interpolation formula}

In this section,  we shall focus our attention on the proof of
Theorem \ref{maint}. To this aim, we shall first introduce several
elementary properties  of the
$c$-divided difference operators.
\begin{lem}\label{lemLeibnitz}  $c$-divided difference operators satisfy the following Leibnitz rules:
\begin{align*}\label{Le}
f(x_1)g(x_1){\cp_1}&=f(x_1)g(x_1){\cp_1}+f(x_1){\cp_1}g(x_2)\\
f(x_1)g(x_1){\cp_1}\cdots{\cp_n}&=\sum_{k=0}^n f(x_1){\cp_1}\cdots{\cp_k}g(x_{k+1}){\cp_{k+1}}
\cdots{\cp_n}.
\end{align*}
\end{lem}

\begin{lem}\label{alem} If $f(x_1,x_2,\ldots,x_n)$ is a symmetric function of $x_i,x_{i+1}$, then
$$
 f(x_1,x_2,\ldots,x_n) {\cp_i}=0.
$$
\end{lem}

\begin{lem} \label{cslem} If $p_n(y)$ is a polynomial of degree $2n$ such that
$$
y^{-n} p_n(y)=(c/y)^{-n}p_{n}(c/y),
$$
then we have
\begin{equation}\label{cslem1}
y_1^{-n} p_n(y_1){\cp_1}\cdots {\cp_m}=\left\{
\begin{array}{ll}
0, &m>n,\\
1, &m=n.
\end{array}
\right.
\end{equation}
\end{lem}
\pf  Clearly, $y_1^{-n}p_n(y_1)$  can be rewritten as $\prod_{i=1}^n(y_1-x_1)(1-c/y_1x_1)$.  When $m=1$ and $n=0$ or $n=1$,  it is easy to verify that \eqref{cslem1} holds.
In view of Lemma  \ref{lemLeibnitz},  we can prove Lemma \ref{cslem} by induction on the length of the operators. \qed

\begin{lem}\label{mlem1}We have
\begin{equation} \label{mlem1e}
\prod_{k=1}^j(y_1-b_k)(1-c/y_1b_k){\cp_1}\cdots{\cp_i}\Big|_{y_k=b_k, 1\leq k \leq i+1}=
\left\{
\begin{array}{ll}
0, &j\neq i;\\
1, &j=i.\\
\end{array}
\right.
\end{equation}
\end{lem}
\pf  For $j\leq i$, Lemma \ref{mlem1} is a direct consequence of Lemma \ref{cslem}. For $j>i$, we have
\begin{align*}
&\prod_{k=1}^j(y_1-b_k)(1-c/y_1b_k){\cp_1}\cdots{\cp_i}\Big|_{y_k=b_k, 1\leq k \leq i+1}\\
&\quad =\prod_{k=2}^j(y_2-b_k)(1-c/y_2b_k){\cp_2}\cdots{\cp_i}\Big|_{y_k=b_k, 2\leq k \leq i+1}\\
& \quad =\cdots = \prod_{k=i+1}^j(y_{i+1}-b_k)(1-c/y_{i+1}b_k)\Big|_{y_k=b_k, i+1\leq k \leq i+1}=0.
\end{align*}
We complete the proof. \qed

\begin{lem} \label{lema} We have
\begin{multline}\label{mlem2}
\prod_{k=1}^j\frac{(y_1-b_k)(1-c/y_1b_k)}{(y_1-a_k)(1-c/y_1a_k)}\prod_{k=1}^{i-1} (y_1-a_k)(1-c/y_1a_k){\cp_1}\cdots{\cp_i}\Big|_{y_k=b_k, 1\leq k \leq i+1}\\
=\left\{
\begin{array}{ll}
0, &j\neq i;\\
\frac{1}{(b_{j+1}-a_j)(1-c/a_jb_{j+1})}, &j=i.\\
\end{array}\right.
\end{multline}

\end{lem}

\pf For $j<i$, we have
\begin{align*}
&\prod_{k=1}^j(y_1-b_k)(1-c/y_1b_k)\prod_{k=j+1}^{i-1} (y_1-a_k)(1-c/y_1a_k){\cp_1}\cdots{\cp_i}\Big|_{y_k=b_k, 1\leq k \leq i+1}\\
&=\sum_{l=0}^i \prod_{k=1}^j(y_1-b_k)(1-c/y_1b_k){\cp_1}\cdots{\cp_l}\Big|_{y_k=b_k, 1\leq k \leq l+1}\\
& \times \prod_{k=l+1}^{i-1} (y_{l+1}-a_k)(1-c/y_{l+1}a_k){\cp_{l+1}}\cdots{\cp_i}\Big|_{y_k=b_k, l+1\leq k \leq i+1}.
\end{align*}

From Lemma \ref{cslem} and Lemma \ref{mlem1}, either the first product inside the sum or the second one vanishes, so does the sum.

For  $j\geq i$,  we have
\begin{align*}
&\frac{\prod_{k=1}^j(y_1-b_k)(1-c/y_1b_k)}{\prod_{k=i}^j(y_1-a_k)(1-c/y_1a_k)}
{\cp_1}\cdots{\cp_i}\Big|_{y_k=b_k, 1\leq k \leq i+1}\\
&=(y_1-b_1)(1-c/yb_1){\cp_1}\Big|_{y_1=b_1}\frac{\prod_{k=2}^j(y_2-b_k)(1-c/y_2b_k)}{\prod_{k=i}^j(y_2-a_k)(1-c/y_2a_k)}
{\cp_2}\cdots{\cp_i}\Big|_{y_k=b_k, 2\leq k \leq i+1}\\
&=\frac{\prod_{k=2}^j(y_2-b_k)(1-c/y_2b_k)}{\prod_{k=i}^j(y_2-a_k)(1-c/y_2a_k)}
{\cp_2}\cdots{\cp_i}\Big|_{y_k=b_k, 2\leq k \leq i+1}\\
&=\cdots=\frac{\prod_{k=i}^j(y_i-b_k)(1-c/y_ib_k)}{\prod_{k=i}^j(y_i-a_k)(1-c/y_ia_k)}
{\cp_i}\Big|_{y_k=b_k, i\leq k \leq i+1}\\
&=\left\{\begin{array}{ll}
0, &j>i,\\
1/(b_{j+1}-a_j)(1-c/a_jb_{j+1}), &j=i.
\end{array}\right.
\end{align*}
We complete the proof.
\qed

\noindent{\it Proof of Theorem \ref{maint}}

Given a polynomial $f(y)$ of degree $2n$ with symmetry $y^{-n}f(y)=(c/y)^{-n}f(c/y)$,  we assume that
\begin{equation}\label{KD}
f(y)=\sum_{j=0}^n C_j \prod_{k=1}^j (y-b_k)(y-c/b_k)\prod_{k=j+1}^{n}(y-a_k)(y-c/a_k).
\end{equation}

Taking $y=b_1$ in \eqref{KD}, one has
$$
f(b_1)=C_0\prod_{k=1}^n (b_1-a_k)(b_1-c/a_k).
$$
Therefore,
$$
C_0=\frac{f(b_1)}{\prod_{k=1}^n (b_1-a_k)(b_1-c/a_k)}.
$$

Setting $y=a_n$ in \eqref{KD} leads to
$$
C_n=\frac{f(a_n)}{\prod_{k=1}^n (a_n-b_k)(a_n-c/b_k)}.
$$

Let $g(y)=f(y)/\prod_{k=1}^n(y-a_k)(y-c/a_k)$. Rewrite \eqref{KD} as
\begin{multline*}
g(y)=g(b_1)+\sum_{j=1}^{n-1} C_j \frac{\prod_{k=1}^j (y-b_k)(1-c/yb_k)}{\prod_{k=1}^j (y-a_k)(1-c/ya_k)}\\
+\frac{f(a_n)}{\prod_{i=1}^n(a_n-b_i)(a_n-c/b_i)}\prod_{i=1}^n\frac{(y-b_k)(1-c/yb_k)}{(y-a_k)(1-c/ya_k)}.
\end{multline*}

Multiplying both sides by $\prod_{k=1}^{i-1}(y-a_k)(1-c/ya_k)$, then applying the operator ${\cp_1}\cdots{\cp_i}$, one has
\begin{align*}
&g(y_1)\prod_{k=1}^{i-1}(y_1-a_k)(1-c/y_1a_k){\cp_1}\cdots{\cp_i}\Big|_{y_k=b_k, 1\leq k \leq i+1}\\
&=\sum_{j=1}^{i-1}C_j\prod_{k=1}^j(y_1-b_k)(1-c/y_1b_k)\prod_{k=j+1}^{i-1}(y_1-a_k)(1-c/y_1a_k)
{\cp_1}\cdots{\cp_i}\Big|_{y_k=b_k, 1\leq k \leq i+1}\\
&+\sum_{j=i}^{n-1}C_j \frac{\prod_{k=1}^j(y_1-b_k)(1-c/y_1b_k)}{\prod_{k=i}^j(y_1-a_k)(1-c/y_1a_k)}{\cp_1}\cdots{\cp_i}\Big|_{y_k=b_k, 1\leq k \leq i+1}.
\end{align*}

By Lemma \ref{lema}, we have
$$
g(y_1)\prod_{k=1}^{i-1}(y_1-a_k)(1-c/y_1a_k){\cp_1}\cdots{\cp_i}\Big|_{y_k=b_k,
1\leq k \leq i+1} =\frac{C_i}{(b_{i+1}-a_i)(1-c/b_{i+1}a_i)}.
$$
Thus
\begin{align*}
C_i&=g(y_1)\prod_{k=1}^{i-1}(y_1-a_k)(1-c/y_1a_k){\cp_1}\cdots{\cp_i}\Big|_{y_k=b_k, 1\leq k \leq i+1}(b_{i+1}-a_i)(1-c/b_{i+1}a_i)\\
&=\frac{f(b_1)b_1^{1-i}}{\prod_{k=i}^n(b_1-a_k)(b_1-c/a_k)}{\cp_1}\cdots{\cp_i}(b_{i+1}-a_i)(1-c/b_{i+1}a_i).
\end{align*}

Replacing $a_i$ by $a_{n-i+1}$, we complete the proof. \qed

\section{Jackson's  $_8\phi_7$   terminating summation formula}

Letting $x_i=uq^{i-1}$ for $1 \leq i \leq n$ and $x_{n+1}=y$, we shall show that \eqref{Kara} in this case is equivalent to Corollary \ref{coro1}.
In other words,  we shall give Jackson's  $_8\phi_7$   terminating summation formula  a determinantal interpretation.

Our proofs in this section involve the following well-known symmetric functions.  Given two sets of variables $X$ and $Y$, the $i$-th  supersymmetric  complete function $h_i(X-Y)$ is defined by
\begin{equation}\label{ss}
 h_{i}(X-Y)=[t^i] \frac{\prod_{y \in Y}(1-yt)}{\prod_{x \in X}(1-xt)}=\sum_{i=0}^n (-1)^{i}e_i(Y)h_{n-i}(X),
\end{equation}
where $[t^i]f(t)$ means the coefficient of $t^i$ in $f(t)$,  $e_i(X)$ and $h_i(Y)$ are $i$-th elementary symmetric function and $i$-th complete symmetric function, respectively.

Expanding the determinant of \eqref{Kara}  along the last row in the case $x_i=uq^{i-1}$ for $1 \leq i \leq n$ and $x_{n+1}=y$ , we have
\begin{eqnarray}\label{87}
\lefteqn{\prod_{i=1}^n(uq^{i-1}-y)(c-uq^{i-1}y)b^{n+1 \choose 2}q^{-(n+1)n(n-1)/3}}\nonumber \\
&&\times \prod_{1 \leq i < j \leq n}(uq^{i-1}-uq^{j-1})(c-u^2q^{i+j-2})\prod_{i=1}^{n+1}(a/b,cq^{2n+2-2i}/ab;q)_{i-1} \nonumber\\
&& =\sum_{k=1}^{n+1}C_{n,k}P_{n-k+1}(y,aq^{k-n})P_{n-k+1}(y,c/a)P_{k-1}(y,bq^{1-n})P_{k-1}(y,cq^{n-k+1}/b), \nonumber\\
\end{eqnarray}
where $C_{n,k}$ is the cofactor of the entry
$$
P_{n-k+1}(y,aq^{k-n})P_{n-k+1}(y,c/a)P_{k-1}(y,bq^{1-n})P_{k-1}(y,cq^{n-k+1}/b).
$$

It is easy to verify that  $C_{n,k}$ can be rewritten in terms of the supersymmetric complete
functions:
\begin{equation}\label{cuab}
C_{n,k}=\prod_{1 \leq i < j \leq n}(uq^{i-1}-uq^{j-1})\det(h_{2n-i+1}(U-Y_{j,k})),
\end{equation}
where the set
$$Y_{j,k}=
\left\{
\begin{array}{ll}
\{a, c/a, \ldots, aq^{j-n}, cq^{n-j}/a, bq^{1-n}, cq^{n-1}/b, \ldots, bq^{j-n-1}, cq^{n-j+1}/b\}, &{\rm if} \, 1\leq j<k,\\
\{a, c/a, \ldots, aq^{j-n+1}, cq^{n-j-1}/a, bq^{1-n},  cq^{n-1}/b, \ldots, bq^{j-n}, cq^{n-j}/b\}, &{\rm if} \, k\leq j\leq n.
\end{array} \right.
$$

For convenience,  we denote by $F_{n,k}(U,A,B)$  the determinant in \eqref{cuab}.  Now, in order to prove  Corollary \ref{coro1}, we are left to evaluate these determinants
$F_{n,k}(U,A,B)$ for $1\leq k\leq n+1$.
\begin{thm}\label{8phi7l}  For $1\leq k \leq n+1$, we have
\begin{align}
 F_{n,k}(U,A,B)  &= {n \brack k-1}
q^{-\frac{n(n-1)(2n-1)}{6}}b^{\binom{n}{2}}\prod_{1\leq i<j \leq n}(c-u^2q^{i+j-2})\nonumber\\
 &\times P_{n-k+1}(b,uq^{n-1})P_{k-1}(a,uq^{n-k+1})P_{n-k+1}(u,c/b)
 P_{k-1}(c/a,u)\nonumber\\
&\times \prod_{i=1}^{n-1}(1-aq^{i-1}/b)^{n-i}\prod_{ j=0 \atop j \neq
n-k+1}^{n-1} \prod_{ i=j+1 \atop i \neq n-k+1}
^{n}(1-cq^{i+j-1}/ab).
\end{align}
\end{thm}

To make the proof clear,  we shall first give two lemmas.

\begin{lem}\label{Las} \cite{Lascoux}
Let $\{j_1,j_2,\ldots,j_n\}$ be a sequence of integers, and let
$X_1$, $\ldots$, $X_n$ and $Y_1$, $\ldots$, $Y_n$ be sets of
variables. The following relation holds
\begin{equation*}
\det\Big(h_{j_k+k-l}(X_k-Y_k)\Big)_{k,l=1}^n
=\det\Big(h_{j_k+k-l}(X_k-Y_k-D_{n-k})\Big)_{k,l=1}^n,
\end{equation*}
where $D_0$, $D_1$, $\ldots$, $D_{n-1}$  are sets of indeterminates
such that the cardinality of $D_i$ is equal to or less than $i$.
\end{lem}

The second lemma is a special case of Theorem \ref{8phi7l}.
\begin{lem} \label{lemm1} For $1\leq k \leq n+1$, we have
\begin{multline}\label{F_{n,k}}
\det(e_{2n-i+1}(Y_{j,k}))_{i,j=1}^n\\
={n \brack k-1}c^{n+1 \choose 2}q^{{n-k+1 \choose
2}-\frac{n(n-1)(2n-1)}{6}}
\prod_{i=1}^{n-1}(b-aq^{i-1})^{n-i}\prod_{ j=0 \atop j \neq
n-k+1}^{n-1} \prod_{ i=j+1 \atop i \neq n-k+1}
^{n}(1-cq^{i+j-1}/ab),
\end{multline}
where the sets of $2n$ variables $Y_{j,k}$, $1\leq j\leq n$, are defined as above.
\end{lem}

\pf  Obviously, $e_{2n}(Y_{j,k})=c^n$. We shall use induction on $n$.
 When $n=1$,  we have
$$
\det(e_2(a,c/a))=\det(e_2(b,c/b))=c.
$$
 Thus \eqref{F_{n,k}} is true for $n=1$.  Assume that \eqref{F_{n,k}} holds for $1 \leq m \leq n-1$, where $n \geq 2$.
 We now proceed to check that \eqref{F_{n,k}} is true for $m=n$.

When $k=1$,  the substraction of two successive columns of the determinant gives
$$
e_i(Y_{j,1})-e_i(Y_{j-1,1})=q^{j-n}(b-a)(1-cq^{2n-2j}/ab)e_{i-1}(Y'_{j-1,1}),
$$
where $Y'_{j-1,1}=Y_{j,1}\setminus\{bq^{j-n},cq^{n-j}/b\}$.

It is easy to verify that
$$
\det(e_{2n-i+1}(Y_{j,1}))_{i,j=1}^n=c^{n+1 \choose 2}q^{{n \choose 2}-\frac{n(n-1)(2n-1)}{6}}\prod_{i=1}^{n-1}(b-aq^{i-1})^{n-i}\prod_{j=0}^{n-1}\prod_{i=j+1}^n (1-cq^{i+j-1}/ab),
$$
which is equal to the right side of \eqref{F_{n,k}}. The case $k=n+1$ is similar.

We now consider the case $2 \leq k \leq n$.  According to
\begin{multline*}
e_i(Y_{k,k})-e_i(Y_{k-1,k})\\
=q^{k-n}(b-a)(1-cq^{2n-2k}/ab)e_{i-1}(Y'_{k-1,k-1})+
q^{k-n-1}(b-a)(1-cq^{2n-2k+2}/ab)e_{i-1}(Y'_{k-1,k}),
\end{multline*}
we have
\begin{multline}
\det(e_{2n-i+1}(Y_{j,k})) = c^nq^{-{n \choose 2}}(b-a)^{n-1}\prod_{i=0}^{n-1}(1-cq^{2i}/ab) \\
 \times \Bigg(\frac{\det(e_{2n-i-1}(Y'_{j,k-1}))q^{n-k}}{1-cq^{2n-2k}/ab}+\frac{\det(e_{2n-i-1}(Y'_{j,k}))q^{n-k+1}}{1-cq^{2n-2k+2}/ab}\Bigg),
\end{multline}
where for $i=k-1$ or $i=k$, we have
$$
Y'_{j,i}=\left\{
\begin{array}{ll}
Y_{j,k}\setminus\{aq^{j-n},cq^{n-j}/a\}, &1\leq j< i,\\
Y_{j+1,k}\setminus\{bq^{j+1-n},cq^{n-j-1}/b\},&i \leq j \leq n-1.
\end{array}
\right.
$$

Our induction hypothesis implies that
\begin{multline*}
\det(e_{2n-i-1}(Y'_{j,k})={n-1 \brack k-1}c^{n \choose 2}q^{{n-k \choose 2}-\frac{(n-1)(n-2)(2n-3)}{6}}\\
\times \prod_{i=1}^{n-2}(bq^{-1}-aq^{i-1})^{n-i-1}\prod_{
j=0 \atop j \neq n-k
}^{n-2}\prod_{i=j+1 \atop i \neq n-k} ^{n-1}(1-cq^{i+j}/ab),
\end{multline*}
and
\begin{multline*}
\det(e_{2n-i-1}(Y'_{j,k-1})={n-1 \brack k-2}c^{n \choose 2}q^{{n-k+1 \choose 2}-\frac{(n-1)(n-2)(2n-3)}{6}}\\
\times \prod_{i=1}^{n-2}(bq^{-1}-aq^{i-1})^{n-i-1}\prod_{
j=0 \atop j \neq n-k+1}^{n-2}
\prod_{
i=j+1 \atop i \neq n-k+1} ^{n-1}(1-cq^{i+j}/ab).
\end{multline*}

Therefore,
\begin{align*}
&\det(e_{2n-i+1}(Y_{j,k}))_{i,j=1}^n\\
&\quad =  c^{n+1 \choose 2}q^{{n-k+1 \choose 2}-\frac{n(n-1)(2n-1)}{6}}\prod_{i=1}^{n-1}(b-aq^{i-1})^{n-i}\prod_{
j=0 \atop j \neq n-k+1}^{n-1}
\prod_{
i=j+1 \atop i \neq n-k+1} ^{n}(1-cq^{i+j-1}/ab) \\
&\quad  \times \Bigg({n-1 \brack k-1}\frac{1-cq^{2n-k}/ab}{1-cq^{2n-2k+1}/ab}+{n-1 \brack k-2}\frac{1-cq^{n-k}/ab}{1-cq^{2n-2k+1}/ab}q^{n-k+1}\Bigg).
\end{align*}
With the aid of the following recurrence
$${n \brack k-1}={n-1 \brack k-1}+q^{n-k+1}{n-1 \brack k-2}=q^{k-1}{n-1 \brack k-1}+{n-1 \brack k-2},$$
we complete the proof.\qed

We are now ready to complete the proof of Theorem \ref{8phi7l}.

\noindent{\it Proof of Theorem \ref{8phi7l}.}

View $F_{n,k}(U,A,B)$ as a polynomial in $u$ of degree $n^2+n$ with coefficients expressed in terms of the other variables.   Applying Lemma \ref{Las},  we first prove that  $F_{n,k}(U,A,B)$ has  $2n $ roots:
$$aq^{1-n}, \ldots, aq^{k-1-n},
cq^{2-k}/a,\ldots,c/a,bq^{k-2n+1},
\ldots,bq^{1-n},c/b,\ldots,cq^{n-k}/b.
$$

Let $u=aq^{i-n}$ in $F_{n,k}(U,A,B)$, where $1 \leq i \leq k-1$.  We  take
$$
D_0=\varnothing, \quad D_1=\{aq\}, \ldots,
D_{i-1}=\{aq,aq^2,\ldots, aq^{i-1}\},
$$
$$
D_{i}=\{aq^{i-n},aq,aq^2,\ldots,aq^{i-1}\}, \ldots,
D_{n-1}=\{aq^{i-n},\ldots,aq^{-1},aq,aq^2,\ldots,aq^{i-1}\},
$$
then apply Lemma \ref{Las}.

Since $e_k(X)=0$ if the cardinality of $X$ is less than $k$,  $F_{n,k}(U,A,B)$ can be transformed
into a determinant whose $(i,j)$-th entry is equal to $0$ if
\[(i,j) \in \{(i,j): 1 \leq j \leq k-1 \hskip 0.2cm \mbox{and}\hskip 0.2cm 1 \leq i \leq n-k+2, \hskip 0.2cm \mbox{or} \hskip 0.2cm k \leq j \leq n-1 \hskip 0.2cm \mbox{and}\hskip 0.2cm 1 \leq i \leq n-j\}.\]
Thus $F_{n,k}(U,A,B)\mid_{u=aq^{i-n}}=0$ for $1\leq i \leq k-1$.

The case $u=cq^{1-i}/a$, where  $1 \leq i \leq k-1$, is similar to the above if we take
$$
D_0=\varnothing, \quad D_1=\{cq^{-1}/a\}, \ldots,
D_{i-1}=\{cq^{1-i}/a,\ldots, cq^{-1}/a\},
$$
$$
D_{i}=\{cq^{1-i}/a,\ldots,cq^{-1}/a,cq^{n-i}/a\}, \ldots,
D_{n-1}=\{cq^{1-i}/a,\ldots,cq^{-1}/a,cq/a,\ldots,cq^{n-i}/a\}.
$$

For the cases $u=bq^{2-n-i}$  and $u=cq^{i-1}/b$,  where $1 \leq i \leq n-k+1$, we take
$$
D_0=\varnothing, \quad D_1=\{bq^{-n}\}, \ldots,
D_{i-1}=\{bq^{2-i-n},\ldots, bq^{-n}\},
$$
$$
D_{i}=\{bq^{2-i-n},\ldots, bq^{-n},bq^{1-i}\}, \ldots,
D_{n-1}=\{bq^{2-i-n},\ldots, bq^{-n},bq^{2-n},\ldots,bq^{1-i}\}
$$
and
$$
D_0=\varnothing, \quad D_1=\{cq^n/b\}, \ldots,
D_{i-1}=\{cq^n/b,\ldots, cq^{n+i-2}/b\},
$$
$$
D_{i}=\{cq^{i-1}/b,cq^n/b,\ldots, cq^{n+i-2}/b\}, \ldots,
D_{n-1}=\{cq^{i-1}/b,\ldots,cq^{n-2}/b,cq^n/b,\ldots,
cq^{n+i-2}/b\},
$$
respectively.

In view of  Lemma \ref{Las}, $F_{n,k}(U,A,B)$  in both cases becomes a
determinant whose $(i,j)$-th entry is equal to $0$ if
\[(i,j) \in \{(i,j): 2 \leq j \leq k-1 \hskip 0.2cm \mbox{and}\hskip 0.2cm 1 \leq i \leq j-1, \hskip 0.2cm \mbox{or} \hskip 0.2cm k \leq j \leq n \hskip 0.2cm \mbox{and}\hskip 0.2cm 1 \leq i \leq k\}.\]
Therefore
$
F_{n,k}(U,A,B)\mid_{u=bq^{2-n-i}}=F_{n,k}(U,A,B)\mid_{u=cq^{i-1}/b}=0
$
for $1 \leq i \leq n-k+1$.

Secondly,  we show that $\prod_{1\leq i<j \leq n}(c-u^2q^{i+j-2})$ is a factor of $F_{n,k}(U,A,B)$.  This is because the determinant vanishes if we set $x_i=c/x_j$ in \eqref{Kara} for each pair $i,j$ where $1\leq i <j \leq n$.

Based on the above, we may assume that
\begin{multline*}
F_{n,k}(U,A,B)=C\times \prod_{1\leq i<j \leq n}(c-u^2q^{i+j-2})\\
\times P_{n-k+1}(b,uq^{n-1})P_{k-1}(a,uq^{n-k+1})P_{n-k+1}(u,c/b)P_{k-1}(c/a,u).
\end{multline*}
To complete the proof, we need to determine $C$. Setting $u=0$, then applying \eqref{F_{n,k}}, we have
\begin{multline*}
C = (-1)^{n-k+1}\frac{\det(e_{2n-i+1}(-Y_{j,k}))_{i,j=1}^n}{c^{n+1
\choose 2}q^{n-k+1 \choose 2}}={n \brack k-1} q^{-\frac{n(n-1)(2n-1)}{6}}\\
\times \prod_{i=1}^{n-1}(b-aq^{i-1})^{n-i}\prod_{ j=0 \atop j \neq
n-k+1}^{n-1} \prod_{ i=j+1 \atop i \neq n-k+1}
^{n}(1-cq^{i+j-1}/ab),
\end{multline*}
as desired. \qed

Putting \eqref{cuab} into \eqref{87},  then divided both sides of
\eqref{87} by
$$
\prod_{1 \leq i < j \leq n}(uq^{i-1}-uq^{j-1})(c-u^2q^{i+j-2})b^{n+1 \choose 2}q^{-(n+1)n(n-1)/3}\prod_{i=1}^{n+1}(a/b,cq^{2n+2-2i}/ab;q)_{i-1},
$$
we complete the proof of Corollary \ref{coro1}.

\vskip 5mm \noindent{\bf Acknowledgments.} This work was supported
by the PCSIRT Project of the Ministry of Education, and the National
Science Foundation of China.

\end{document}